\documentclass[10pt]{amsart}
\usepackage{latexsym}
\usepackage{amsfonts}

\usepackage{amsmath,amsthm,amssymb}

\title[On the absence of McShane-type identities for the outer space]
{On the absence of McShane-type identities for the outer space}

\author[I.~Kapovich]{Ilya Kapovich}

\address{\tt Department of Mathematics, University of Illinois at
  Urbana-Champaign, 1409 West Green Street, Urbana, IL 61801, USA
  \newline http://www.math.uiuc.edu/\~{}kapovich/} \email{\tt
  kapovich@math.uiuc.edu}

\curraddr{\tt Department of Mathematics, Johann Wolgang Goethe
  University, Robert Mayer-Strasse 6-8,
60325 Frankfurt (Main) , Germany}
\email{\tt kapovich@math.uni-frankfurt.de}

\author[I.~Rivin]{Igor Rivin} \address{\tt Department of Mathematics,
  Temple University, Philadelphia, PA 19122, USA} \email{\tt rivin@euclid.math.temple.edu}

\newtheorem{theor}{Theorem}

\newtheorem{thm}{Theorem}[section] 
 
\newtheorem{prop}[thm]{Proposition} \theoremstyle{definition}
\newtheorem{defn}[thm]{Definition}

\newtheorem{conv}[thm]{Convention} \newtheorem{rem}[thm]{Remark}

\begin{document}

\begin{abstract}
A remarkable result of McShane states that for a punctured torus
with a complete finite volume hyperbolic metric we have
\[
\sum_{\gamma} \frac{1}{e^{\ell(\gamma)}+1}=\frac{1}{2}
\]
where $\gamma$ varies over the homotopy classes of essential
simple closed curves and $\ell(\gamma)$ is the length of the
geodesic representative of $\gamma$.

We prove that there is no reasonable analogue of McShane's
identity for the Culler-Vogtmann outer space of a free group.
\end{abstract}


\thanks{The first author was supported by the NSF
  grant DMS\#0404991 and by the Humboldt Foundation Fellowship}

\maketitle

\section{Introduction}\label{intro}

Let $T$ be the one-punctured torus and let $\rho$ be a complete
finite-volume hyperbolic structure on $T$. Let $\mathcal S$ be the
set of all free homotopy classes of essential simple closed curves
in $T$ that are not homotopic to the puncture. Denote
\[
E(\rho):=\sum_{\gamma \in \mathcal S}
\frac{1}{e^{\ell_\rho(\gamma)}+1},
\]
where $\ell_\rho(\gamma)$ is the smallest $\rho$-length among all
curves representing $\gamma$. Thus $E$ can be regarded as a
function on the Teichm\"uller space of $T$. A remarkable result of
McShane~\cite{Mc91} shows that this function is constant and that
\begin{gather*}
E(\rho)=\frac{1}{2}\tag{$\ast$}
\end{gather*}
for every $\rho$. We refer to $(\ast)$ as \emph{McShane's
identity} for $T$. Since the thesis of McShane~\cite{Mc91}, other
proofs of McShane's identity for the punctured torus have been
produced (particularly, see the work of Bowditch~\cite{Bow96}) and
McShane's identity has been generalized to other hyperbolic
surfaces and other contexts~\cite{Bow97,Mc98,AMS1,AMS2,TWZ,TWZ1}. Note that
if $\psi$ is an element of the mapping class group of $T$ then
$\psi$ permutes the elements of $\mathcal S$ and hence, clearly,
$E(\rho)=E(\psi\rho)$. Thus $E$ obviously factors through to a
function on the moduli space of $T$ and $(\ast)$ says that this
function is identically equal to $1/2$.

Let $F_k=F(a_1,\dots, a_k)$ be a free group of rank $k\ge 2$ with a free
basis $A=\{a_1\dots, a_k\}$.  For $F_k$ the best analogue of the
Teichm\"uller space is the so-called Culler-Vogtmann \emph{outer space}
$CV(F_k)$. Instead of actions on the hyperbolic plane the elements of
the outer space are represented by minimal discrete isometric actions of
$F_k$ on $\mathbb R$-trees.  Equivalently, one can think about a point
of the outer space as being represented by a \emph{marked volume-one
  metric graph structure} on $F_k$, that is, an isomorphism $\phi:
F_k\to \pi_1(\Gamma, p)$, where $\Gamma$ is a finite graph without
degree-one and degree-two vertices, equipped with a \emph{metric
  structure} $\mathcal L$ that assigns to each non-oriented edge of
$\Gamma$ a positive number called the \emph{length} of this edge.  The
\emph{volume} of a metric structure on $\Gamma$ is the sum of the
lengths of all non-oriented edges of $\Gamma$. As we noted, the metric
structures that appear in the description of the points of the outer
space, given above, are required to have volume equal to one. If
$(\phi:F_k\to \pi_1(\Gamma, p), \mathcal L)$ represents a point of the
outer space, the metric structure $\mathcal L$ naturally lifts to the
universal cover $\widetilde \Gamma$, turning $\widetilde \Gamma$ into an
$\mathbb R$-tree $X$. The group $F_k$ acts on this $\mathbb R$-tree $X$
via $\phi$ by isometries minimally and discretely with the quotient
being equal to $\Gamma$. Similarly to the marked length spectrum in the
Teichm\"uller space context, a marked metric graph structure
$(\phi:F_k\to \pi_1(\Gamma, p), \mathcal L)$ defines a \emph{hyperbolic
  length function} $\ell: \mathcal C_k\to\mathbb R$ where $\mathcal C_k$
is the set of all nontrivial conjugacy classes in $F_k$.  If $g\in F_k$,
then $\ell([g])$ is the translation length of $g$ considered as the
isometry of the $\mathbb R$-tree $X$ described above. Alternatively, we
can think about $\ell([g])$ as follows: $\ell([g])$ is the $\mathcal
L$-length of the shortest free homotopy representative of the curve
$\phi(g)$ in $\Gamma$, that is, the $\mathcal L$-length of the
"cyclically reduced" form of $\phi(g)$ in $\Gamma$. Two volume-one
metric graph structures on $F_k$ represent the same point of $CV(F_k)$
if and only if their corresponding hyperbolic length functions are
equal, or, equivalently, if the corresponding $\mathbb R$-trees are
$F_k$-equivariantly isometric.

It is natural to ask if there is an analogue of McShane's identity
in the outer space context. The (right) action of $\psi\in
Out(F_k)$ on $CV(F_k)$ takes a hyperbolic length function $\ell$
to $\ell\circ \psi$, that is, $\psi$ simply permutes the domain
$\mathcal C_k$ of $\ell$. Therefore the real question, as in the
Teichm\"uller space case, is if there is an analogue of McShane's
identity for the \emph{moduli space} $\mathcal
M_k=CV(F_k)/Out(F_k)$. The elements of $\mathcal M_k$ are
represented by unmarked finite connected volume-one metric graphs
$(\Gamma, \mathcal L)$ without degree-one and degree-two vertices
and with $\pi_1(\Gamma)\simeq F_k$.

To simplify the picture, and also since our results will be negative, we
will consider a subset $\Delta_k$ of $CV(F_k)$ consisting of all
volume-one metric structures on the wedge $W_k$ of $k$ circles wedged at
a base-vertex $v_0$. We orient the circles and label them by $a_1,\dots,
a_k$.  This gives us an identification $\pi_1(W_k,v_0)=F(a_1,\dots,
a_k)$ of $F_k=F(a_1,\dots, a_k)$ with $\pi_1(W_k,v_0)$, so that indeed
$\Delta_k\subseteq CV(F_k)$.

A volume-one metric structure $\mathcal L$ on $W_k$ is a $k$-tuple
$(\mathcal L(a_1), \dots, \mathcal L(a_k))$ of positive numbers
with $\sum_{i=1}^k\mathcal L(a_i)=1$. Thus $\Delta_k$ has the
natural structure of an open $(k-1)$-dimensional simplex in
$\mathbb R^k$. As in the general outer space context, every
$\mathcal L\in \Delta_k$ defines a hyperbolic length-function
$\ell_\mathcal L:\mathcal C_k\to \mathbb R$, where for $g\in F_k$
$\ell_\mathcal L([g])$ is the $\mathcal L$-length of the
cyclically reduced form of $g$ in $F_k=F(a_1,\dots, a_k)$. The
open simplex $\Delta_k$ has a distinguished point $\mathcal
L_\ast:=(\frac{1}{k},\dots, \frac{1}{k})$. Note that for every
$[g]\in \mathcal C_k$ we have $\ell_{\mathcal
L_\ast}([g])=||g||/k$, where $||g||$ is the cyclically reduced
length of $g$ in $F_k=F(a_1,\dots, a_k)$.

There is no perfect  analogue for the notion of a simple
closed curve in the free group context. The closest such analogue
is given by \emph{primitive elements}, that is, elements belonging
to some free basis of $F_k$. Let $\mathcal P_k$ be denote the set
of conjugacy classes of primitive elements of $F_k$. We will
consider two versions of possible generalizations of McShane's
identity for free groups: the first involving all conjugacy
classes in $F_k$ and the second involving the conjugacy classes of
primitive elements of $F_k$. We will see that, under some
reasonable assumptions, there are no analogues of McShane's
identity in either context.

\begin{defn}[McShane-type functions on $\Delta_k$]
Let $f:(0,\infty)\to (0,\infty)$ be a monotone non-increasing
function and let $k\ge 2$. Define

\[
C_f:\Delta_k\to (0, \infty], \qquad C_f(\mathcal L)=\sum_{w\in
\mathcal C_k} f(\ell_\mathcal L(w))\quad \text{ where }\mathcal
L\in \Delta_k
\]
and
\[
P_f:\Delta_k\to (0, \infty], \qquad P_f(\mathcal L)=\sum_{w\in
\mathcal P_k} f(\ell_\mathcal L(w))\quad \text{ where }\mathcal
L\in \Delta_k.
\]
\end{defn}
Obviously, $0<P_f<C_f\le \infty$ on $\Delta_k$.

Motivated by McShane's identity, it is interesting to ask if there
exist functions $f$ such that either $C_f$ or $P_f$ is constant on
$\Delta_k$. To make the question meaningful we need to require
$P_f$ (or, correspondingly, $C_f$) be finite at some point
$\mathcal L\in \Delta_k$. Thus it is necessary to assume that
$\lim_{x\to\infty} f(x)=0$ and that this convergence to zero is
sufficiently fast.

We establish the following negative results regarding the
existence of analogues of McShane's identity in the outer space
context:

\begin{theor}\label{A}
Let $k\ge 2$ be an integer and let $F=F(a_1,\dots, a_k)$. Let $f:
(0,\infty)\to (0,\infty)$ be a monotone non-increasing function
such that:
\begin{enumerate}
\item
 \[
\limsup_{x\to \infty} f(x)^{1/x}< \frac{1}{(2k-1)^k}.
\]
\item
\[
\liminf_{x\to \infty} f(x)^{1/x}> 0.
\]
\end{enumerate}

Then:

\begin{enumerate}
\item[(a)] We have $P_f\le C_g<\infty$ on some neighborhood $U$ on
$\mathcal L_\ast$ in $\Delta_k$ (moreover, only the assumption (1)
on $f$ above is required for this conclusion).

\item[(b)] We have $C_f\ne const$ on $\Delta_k$.

\item[(c)] If $k\ge 3$ then $P_f\ne const$ on $\Delta_k$.
\end{enumerate}

\end{theor}

The assumptions on $f(x)$ in Theorem~\ref{A} require $f(x)$ to
decay both at least and at most exponentially fast; condition (1)
assures that the value of $C_f$ is finite near $\mathcal L_\ast$.
The idea of the proof of parts (b) and (c) of Theorem~\ref{A} uses
the notion of \emph{volume entropy} for a metric structure
$\mathcal L$ on $W_k$ (see  ~\cite{KN,Riv99,L}). Roughly speaking, there are points
$\mathcal L$ near the the boundary of $\Delta_k$ where the
exponential growth rate, as $R\to\infty$, of the number of
conjugacy classes with $\ell_\mathcal L$-length at most $R$ is
bigger than the exponential rate of decay of the function $f$.
This forces $C_f$ to be equal to $\infty$ at $\mathcal L$.

For $k=2$ the set of conjugacy classes of primitive elements has
quadratic rather than exponential growth. Therefore we modify the
assumptions on $f(x)$ accordingly and obtain a somewhat stronger
conclusion then in part (c) of Theorem~\ref{A}. For $k=2$ the open
$1$-dimensional simplex $\Delta_2\subseteq \mathbb R^2$ consists
of all pairs $\mathcal L_t:=(t,1-t)$ where $t\in (0,1)$. Therefore
we may identify $\Delta_2$ with $(0,1)$ and define
$P_f(t):=P_f(\mathcal L_t)$. With this convention we prove:

\begin{theor}\label{B}
Let $k=2$ and $F=F(a,b)$. Let $f: (0,\infty)\to (0,\infty)$ be a
monotone non-increasing function such that:
\begin{enumerate}
\item We have $f''(x)>0$ for every $x>0$.

\item There is some $\epsilon>0$ such that  $\lim_{x\to\infty}
x^{3+\epsilon}f(x) =0$.
\end{enumerate}
Then the following hold:

\begin{enumerate}
\item[(a)] We have $0<P_f(t)<\infty$ for every $t\in (0,1)$.

\item[(b)] The function $P_f(t)$ is strictly convex on $(0,1)$ and
achieves a unique minimum at $t_0=1/2$. In particular, $P_f(t)$ is
not a constant locally near $t_0=1/2$ and thus $P_f\ne const$ on
$(0,1)$.
\end{enumerate}

\end{theor}
The proof of Theorem~\ref{B} uses convexity considerations as well
as some results about the explicit structure of primitive elements
in $F(a,b)$~\cite{CMZ,OZ}.

Finally, we combine the volume entropy and the convexity ideas to
obtain:

\begin{theor}\label{C}
Let $k\ge 2$ and let $f:(0,\infty)\to (0,\infty)$ be a monotone
decreasing function such that the following hold:
\begin{enumerate}
\item The function $f(x)$ is strictly convex on $(0,\infty)$.

\item
\[ \limsup_{x\to \infty} f(x)^{1/x}< \frac{1}{(2k-1)^k}.
\]
\end{enumerate}

Then there exists a convex neighborhood $U$ of $\mathcal L_\ast$
in $\Delta_k$ such that $0<P_f< C_f<\infty$ on $U$ and both $C_f$
and $P_f$ are strictly convex on $U$. In particular, $C_f\ne
const$ on $U$ and $P_f\ne const$ on $U$.
\end{theor}

The authors are grateful to Paul Schupp for useful conversations.

\section{Volume entropy}

In this section we will prove Theorem~\ref{A}, which is obtained
as a combination of Theorem~\ref{thm:conj} and
Theorem~\ref{thm:prim1} below.

\begin{conv}
For the remainder of this section let $k\ge 2$ be an integer and
$F_k=F(a_1,\dots, a_k)$ be free of rank $k$ with a free basis
$A=\{a_1,\dots, a_k\}$. We identify $F_k$ with $\pi_1(W_k,v_0)$,
as explained in the introduction. For $g\in F_k$ we denote by
$|g|$ the freely reduced length of $g$ with respect to $A$ and we
denote by $||g||$ the cyclically reduced length of $g$ with
respect to $A$.

We denote by $CR_k$ the set of all cyclically reduced elements of
$F_k$ with respect to $A$.

Let $\mathcal L$ be a metric graph structure on $W_k$.  For every
$g\in F_k$ there is a unique edge-path in $W_k$ labelled by the
freely reduced form of $g$ with respect to $A$. We denote the
$\mathcal L$-length of that path by $\mathcal L(g)$. As before, we
denote by $\ell_\mathcal L: \mathcal C_k\to\mathbb R$ the
hyperbolic length function corresponding to $\mathcal L$. Thus if
$g\in F_k$ then $\ell_\mathcal L([g])=\mathcal L(u)$ where $u$ is
the cyclically reduced form of $g$ with respect to $A$.
\end{conv}

\begin{defn}[Volume Entropy]\label{defn:ent}
Let $\mathcal L$ be a metric structure on $W_k$.  The \emph{volume
entropy} $h_\mathcal L$ of $\mathcal L$ is defined as
\[
h_\mathcal L=\lim_{R\to\infty}\frac{\log \#\{ g\in F_k: \mathcal
    L(g)\le R \}}{R}.
\]
\end{defn}

It is well-known and easy to see that the limit in the above expression
exists and is finite. We refer the reader to~\cite{KN,Riv99,L} for a
detailed discussion of volume entropy in the context of metric graphs.

\begin{prop}\label{prop:ent}
Let $k\ge 2$ and $\mathcal L$ be as in definition~\ref{defn:ent}.

Then the limits
\[
h'_\mathcal L=\lim_{R\to\infty}\frac{\log \#\{ g\in CR_k: \mathcal
    L(g)\le R \}}{R}.
\]

and

\[
h''_\mathcal L=\lim_{R\to\infty}\frac{\log \#\{ w\in \mathcal C_k:
\ell_{\mathcal
    L}(w)\le R \}}{R}.
\]

exist and
\[
h_\mathcal L=h'_\mathcal L=h''_\mathcal L.
\]
\end{prop}

\begin{proof}

Let $M:=\max\{ |a_i|_\mathcal L: i=1,\dots, k\}$ and $m:=\min\{ |a_i|_\mathcal L: i=1,\dots, k\}$

For each  $g\in F$ there exists a cyclically reduced word $v_g$
such that $|g|=|v_g|$ and such that $g$ and $v_g$ agree except
possibly in the last letter. Then $\big|  \mathcal L(g)-\mathcal
L(v_g)\big |\le M$. Moreover, the function $F_k\to CR_k, g\mapsto
v_g$ is at most $2k$-to-one. Therefore for every integer $R>0$
\[
\#\{ g\in CR_k: \mathcal L(g)\le R \}\le \#\{ g\in F_k: \mathcal
  L(g)\le R \}\le 2k \#\{ g\in CR_k: \mathcal L(g)\le R+M\}
\]
and
\[
\#\{ w\in \mathcal C_k: \ell_{\mathcal L}(w)\le R \}\le \#\{ g\in
CR_k: {\mathcal L}(g)\le R \}\le \frac{R}{m} \#\{ w\in \mathcal
C_k: \ell_{\mathcal L}(w)\le R \}.
\]
This implies the statement of the proposition.
\end{proof}

\begin{thm}\label{thm:conj}
Let $k\ge 2$ be an integer and let $F=F(a_1,\dots, a_k)$. Let $f:
(0,\infty)\to (0,\infty)$ be a monotone non-increasing function
such that:
\begin{enumerate}
\item
 \[
\limsup_{x\to \infty} f(x)^{1/x}< \frac{1}{(2k-1)^k}.
\]
\item
\[
\liminf_{x\to \infty} f(x)^{1/x}> 0.
\]
\end{enumerate}

Then:

\begin{enumerate}
\item[(a)] We have $0<C_f<\infty$ on some neighborhood $U$ on
$\mathcal L_\ast$ in $\Delta_k$ (moreover, only the assumption (1)
on $f$ above is required for this conclusion).

\item[(b)] We have $C_f\ne const$ on $\Delta_k$.
\end{enumerate}

\end{thm}
\begin{proof}
The assumptions on $f(x)$ imply that there exist $N>0$ and
$0<\sigma_1< \sigma_2<\frac{1}{(2k-1)^k}$ such that for every
$x\ge N$
\[
\sigma_1^x \le f(x) \le \sigma_2^x.
\]

Let $\mathcal L_\ast=(\frac{1}{k}, \dots, \frac{1}{k})\in
\Delta_k$. For any $g\in F_k$ we have $\mathcal L_\ast(g)=|g|/k$.
Then an easy  direct computation shows that $h_{\mathcal
L_\ast}=k\log(2k-1)$, so that $e^{h_{\mathcal
L_\ast}}=(2k-1)^k<\frac{1}{\sigma_2}$. Since the volume entropy $h$ is a continuous
function on $\Delta_k$ (see, for example, \cite{KN}), there exist
a neighborhood $U$ of $\mathcal L_\ast$ in $\Delta_k$ and
$0<c<\frac{1}{\sigma_2}$ such that $e^{h_{\mathcal L}}<c$ for
every $\mathcal L\in U$.

Observe now that $C_f<\infty$ on $U$. Let $\mathcal L\in U$ be
arbitrary. There there exist $M>0$ and $c_1$ with
$c<c_1<\frac{1}{\sigma_2}$ such that for every integer $R>0$ we
have
\[
\#\{w\in \mathcal C_k: \ell_{\mathcal L}(w)\le R\}\le M c_1^R
\]

Therefore
\begin{gather*}
C_f(\mathcal L_\ast)=\sum_{w\in \mathcal C_k} f(\ell_{\mathcal
L}(w))=\sum_{i=0}^{\infty}\ \ \sum_{w\in \mathcal C_k,
i<\ell_{\mathcal L}(w)\le i+1} f(\ell_{\mathcal L}(w))\le \\
\sum_{i=0}^{\infty}\  \  \sum_{w\in \mathcal C_k, i<\ell_{\mathcal
L}(w)\le i+1} f(i)\le \sum_{i=0}^{\infty}  M c_1^{i+1} f(i)<\infty
\end{gather*}
where the last inequality holds since $c_1<\frac{1}{\sigma_2}$ and
$f(x) \le \sigma_2^x$ for all $x\ge N$. Thus indeed $C_f(\mathcal
L)<\infty$, so that $C_f<\infty$ on $U$.

For $0<t<\frac{1}{k-1}$ let $\mathcal L_t=(t,t,\dots, t,
1-(k-1)t)\in \Delta_k$. Then, as follows from the proof of
Theorem~B of~\cite{KN} (specifically the proof of Theorem~9.4 on page 25
of~\cite{KN} ),
\[
\lim_{t\to 0} h_{\mathcal L_t}=\infty.
\]
 
Indeed, let $\Gamma$ be the subgraph of $W_k$ consisting of the loops
labelled by $a_1$ and $a_k$. The restriction $\mathcal L_t'$ of
$\mathcal L_t$ to $\Gamma$ is a metric structure on $\Gamma$ of volume
$1-(k-2)t$. Therefore $\frac{1}{1-(k-2)t}\mathcal L_t'$ is a volume-one
metric structure on $\Gamma$ with respect to which the length of $a_1$
goes to $0$ as $t\to 0$. Therefore, as established in the proof of
Theorem~9.4 of~\cite{KN}, 
\[
\lim_{t\to 0} h_{\frac{1}{1-(k-2)t}\mathcal L_t'}=\infty.
\] 
However, 
\[
h_{\mathcal L_t'}=\frac{1}{1-(k-2)t}h_{\frac{1}{1-(k-2)t}\mathcal L_t'}
\]
and therefore
\[
\lim_{t\to 0} h_{\mathcal L_t'}=\infty.
\] 
It is obvious from the definition of volume entropy that
$h_{\mathcal L_t'}\le h_{\mathcal L_t}$ and hence
\[
\lim_{t\to 0} h_{\mathcal L_t}=\infty,
\]
as claimed.

Hence there exists $t_0\in (0,\frac{1}{k-1})$ such that for every
$t\in (0,t_0)$ we have $e^{h_{\mathcal
L_t}}>\frac{1}{\sigma_1}+2$. Let $t\in (0,t_0)$ be arbitrary. We
claim that $C_f(\mathcal L_t)=\infty$.

Since $e^{h_{\mathcal L_t}}>\frac{1}{\sigma_1}+2$, by
Proposition~\ref{prop:ent} there is $R_0>N>0$ such that for every
$R\ge R_0$ we have
\[
\#\{w\in \mathcal C_k:  \mathcal L_t(w)\le R\}\ge
(\frac{1}{\sigma_1}+1)^R.
\]
For every $R\ge R_0$
\begin{gather*}
C_f(\mathcal L_t)=\sum_{w\in \mathcal C_k} f(\mathcal L_t(w)) \ge
\sum_{w\in \mathcal C_k,
  \mathcal L_t(w)\le R} f(\mathcal L_t(w)) \ge \\
\sum_{w\in \mathcal C_k, \mathcal L_t(w)\le R} f(R) \ge
(\frac{1}{\sigma_1}+1)^R f(R) \ge
(\frac{1}{\sigma_1}+1)^R\sigma_1^R =\\
=(1+\sigma_1)^R.
\end{gather*}
Since this is true for every $R\ge R_0$, it follows that
$C_f(\mathcal L_t)=\infty$.

Thus $C_f(\mathcal L_\ast)<\infty$ while $C_f(\mathcal
L_t)=\infty$ for all sufficiently small $t>0$. Therefore $C_f\ne
const$ on $\Delta_k$.
\end{proof}

\begin{thm}\label{thm:prim1}
Let $k\ge 3$ be an integer and let $F=F(a_1,\dots, a_k)$. Let $f:
(0,\infty)\to (0,\infty)$ be as in Theorem~\ref{thm:conj}.

Then:

\begin{enumerate}
\item[(a)] We have $P_f\le C_f<\infty$ on some neighborhood $U$ on
$\mathcal L_\ast$ in $\Delta_k$.

\item[(b)] We have $P_f\ne const$ on $\Delta_k$.
\end{enumerate}

\end{thm}

\begin{proof}
Again, by assumptions on $f(x)$, there exist there exist $N>0$ and
$0<\sigma_1< \sigma_2<\frac{1}{(2k-1)^k}$ such that for every
$x\ge N$
\[
\sigma_1^x \le f(x) \le \sigma_2^x.
\]

By Definition $0\le P_f\le C_f$. By Theorem~\ref{thm:conj} we have
$C_f<\infty$ on some neighborhood $U$ on $\mathcal L_\ast$ in
$\Delta_k$ and hence $P_f\le C_f<\infty$ on $U$.

Put $F_{k-1}:=F(a_1,\dots, a_{k-1})$ so that $F_k=F_{k-1}\ast
\langle a_k\rangle$. For $0<t<\frac{1}{k-2}$ let \[\mathcal
L_t:=(\frac{t}{2},\frac{t}{2},\dots, \frac{t}{2},
\frac{1}{2}-(k-2)\frac{t}{2}, \frac{1}{2})\in \Delta_k\] and
\[\widehat {\mathcal L_t}:=(\frac{t}{2},\frac{t}{2},\dots,
\frac{t}{2}, \frac{1}{2}-(k-2)\frac{t}{2})\in
\frac{1}{2}\Delta_{k-1}.\] Thus $2\widehat {\mathcal L_t}\in
\Delta_{k-1}$ is a volume-one metric structure on $W_{k-1}$. Since
$k\ge 3$, we have $k-1\ge 2$ and hence, exactly as in the proof of
Theorem~\ref{thm:conj}, $\lim_{t\to 0} {h_{2\widehat {\mathcal
L_t}}}=\infty$.

For $R\ge 1$
\[
b_{R,t}:=\#\{g\in F_{k-1}: \widehat{\mathcal L_t}(g)\le R\}.
\]
Then
\[
s_t:=\lim_{R\to\infty} \frac{\log b_{R,t}}{R}=h_{\widehat{\mathcal
L_t}}=2h_{2\widehat{\mathcal L_t}}.
\]
and therefore
\[
\lim_{t\to 0} s_t=\infty.
\]

Hence there exists $0<t_0<\frac{1}{k-2}$ such that for every $t\in
(0,t_0)$ we have $e^{s_t}>\frac{1}{\sigma_1}+2$.

 Fix an arbitrary $t\in (0,t_0)$.
Since $e^{s_t}>\frac{1}{\sigma_1}+2$, by there is $R_0>N>0$ such
that for every $R\ge R_0$ we have
\[
b_{R,t}=\#\{g\in F_{k-1}: \widehat{\mathcal L_t}(g)\le R\}\ge
(\frac{1}{\sigma_1}+1)^R.
\]
Note that for every $g\in F_{k-1}$ the element $ga_k\in F_k$ is
primitive in $F$. Moreover, if $g_1\ne g_2$ are distinct elements
of $F_{k-1}$ then $g_1a_k$ and $g_2a_k$ are not conjugate in
$F_k$. Recall that by definition of $\mathcal L_t$ we have
$\mathcal L_t(a_k)=\frac{1}{2}$. For $R\ge 1$ denote
\[
p_{R,t}:=\#\{w\in \mathcal P_k: \ell_{\mathcal L_t}(w)\le R\}.
\]

Then for every $R\ge R_0+\frac{1}{2}$ we have
\[
p_{R,t}\ge b_{R-\frac{1}{2},t}\ge
(\frac{1}{\sigma_1}+1)^{R-\frac{1}{2}}.
\]
Hence for every $R\ge R_0+\frac{1}{2}$
\begin{gather*}
P_f(\mathcal L_t)=\sum_{w\in \mathcal P_k} f(\mathcal L_t(w)) \ge
\sum_{w\in \mathcal P_k,
  \ell_{\mathcal L_t}(w)\le R} f(\mathcal L_t(w)) \ge \\
\sum_{w\in \mathcal P_k, \ell_{\mathcal L_t}(w)\le R} f(R) \ge
(\frac{1}{\sigma_1}+1)^{R-\frac{1}{2}} f(R) \ge
(\frac{1}{\sigma_1}+1)^{R-\frac{1}{2}}\sigma_1^R =\\
=(1+\sigma_1)^R(\frac{1}{\sigma_1}+1)^{-\frac{1}{2}}.
\end{gather*}
Since this is true for every $R\ge R_0+\frac{1}{2}$, it follows
that $P_f(\mathcal L_t)=\infty$.

Thus $P_f(\mathcal L_\ast)<\infty$ while $P_f(\mathcal
L_t)=\infty$ for all sufficiently small $t>0$. Therefore $P_f\ne
const$ on $\Delta_k$.

\end{proof}

\section{Primitive elements in $F(a,b)$}

In this section we will prove Theorem~\ref{B}.

\begin{conv}
Throughout this section let $F_2=F(a,b)$ be a free group of rank
two.

Let $\alpha:F(a,b)\to\mathbb Z^2$ be the abelianization
homomorphism, that is, $\alpha(a)=(1,0)$ and $\alpha(b)=(0,1)$.
Then $\alpha$ is constant on every conjugacy class and therefore
$\alpha$ defines a map $\beta: \mathcal C_2\to \mathbb Z^2$.
\end{conv}

\begin{defn}[Visible points]
A point $(p,q)\in \mathbb Z^2$ is called \emph{visible} if
$gcd(p,q)=1$. We denote the set of all visible points in $\mathbb
Z^2$ by $V$.
\end{defn}

We will need the following known facts about primitive elements in
$F(a,b)$ (see, for example, \cite{CMZ,OZ}:

\begin{prop}\label{prop:prim}
The following hold:
\begin{enumerate}
\item The restriction of $\beta$ to $\mathcal P_2$ is a bijection
between $\mathcal P_2$ and the set of visible elements $V\subseteq
\mathbb Z^2$.

\item Let $w\in F(a,b)$ be a cyclically reduced primitive element
and let $\alpha(w)=(p,q)\in \mathbb Z^2$.

Then every occurrence of $a$ in $w$ has the same sign (either
$-1,0$ or $1$) as $p$ and every occurrence of $b$ in $w$ has the
same sign (again either $-1,0$ or $1$) as $q$. Thus the total
number of occurrences of $a^{\pm 1}$ in $w$ is equal to $|p|$ and
the total number of occurrences of $b^{\pm 1}$ in $w$ is equal to
$|q|$.
\end{enumerate}

\end{prop}

\begin{defn}[Admissible function]
We say that a function $f: (0,\infty)\to [0,\infty)$ is
\emph{admissible} if it satisfies the following conditions:
\begin{enumerate}
\item We have $f''(x)>0$ for every $x>0$.

\item There is some $\epsilon>0$ such that  $\lim_{x\to\infty}
x^{3+\epsilon}f(x) =0$.
\end{enumerate}

\end{defn}

The second condition means that $f(x)$ converges to zero
asymptotically faster than $\frac{1}{x^{3+\epsilon}}$ as
$x\to\infty$. Note that an admissible function must be strictly
positive and monotone decreasing on $(0,\infty)$.

\begin{thm}\label{thm:prim}
Let $f$ be any admissible function. Then the following hold:

\begin{enumerate}
\item We have $0<P_f(t)<\infty$ for every $t\in (0,1)$.

\item The function $P_f(t)$ is strictly convex on $(0,1)$ and
achieves a unique minimum at $t=1/2$. In particular, $P_f(t)$ is
not a constant locally near $t=1/2$.

\end{enumerate}

\end{thm}
\begin{proof}
For every $(p,q)\in V$ and $t\in (0,1)$ denote
\[
g_{p,q}(t)=f(t|p|+(1-t)|q|)+f(t|q|+(1-t)|p|).
\]

Note that if $(p,q)\in V$ then $gcd(p,q)=1$ and hence $|p|\ne
|q|$. We can therefore partition $V$ as the collection of pairs
$(p,q), (q,p)$ of visible elements and every such pair has a
unique representative where the absolute value of the first
coordinate is bigger than that of the second coordinate.

Let $w\in \mathcal P_2$ be arbitrary and let $(p,q)=\beta(w)$.
Proposition~\ref{prop:prim} and the definition of $\mathcal L_t$
imply that for any $t\in (0,1)$ we have
\[
\ell_{\mathcal L_t}(w)=t|p|+(1-t)|q|.
\]

Let $V':=\{(p,q)\in V: |p|>|q|\}$. Then we have

\begin{gather*}
P_f(t)=\sum_{w\in \mathcal P_2} f(\ell_{\mathcal
L_t}(w))=\sum_{(p,q)\in
V}f(t|p|+(1-t)|q|)=\tag{\dag}\\
\sum_{(p,q)\in V'}f(t|p|+(1-t)|q|)+f(t|q|+(1-t)|p|)=\sum_{(p,q)\in
V'}g_{p,q}(t).
\end{gather*}

Fix some $t\in (0,1)$. We can also represent $P_f(t)$ as
\[
P_f(t)=\sum_{N=1}^{\infty}\ \  \sum_{(p,q)\in V,
\max\{|p|,|q|\}=N}f(t|p|+(1-t)|q|).
\]
 Since $f(x)$ is a monotone
non-increasing function, if $(p,q)\in V,
\max\{|p|,|q|\}=N$, we have
\[
f(t|p|+(1-t)|q|)\le \min\{f(tN), f((1-t)N)\}=f(cN)
\]
where $c=\max\{t, 1-t\}$. For every integer $N\ge 1$ the number of
points $(p,q)\in \mathbb Z^2$ with $|p|\le N, |q|\le N$ is
$(2N+1)^2$.

Therefore
\begin{gather*}
P_f(t)=\sum_{N=1}^{\infty}\ \  \sum_{(p,q)\in V,
\max\{|p|,|q|\}=N}f(t|p|+(1-t)|q|)\le \\
\sum_{N=1}^{\infty}\ \  \sum_{(p,q)\in V, \max\{|p|,|q|\}=N}
f(cN)\le\sum_{N=1}^{\infty} (2N+1)^2 f(cN)<\infty
\end{gather*}
because of condition (2) in the definition of admissibility of
$f(x)$. Thus $0<P_f(t)<\infty$ for every $t\in (0,1)$.

Note that for each $(p,q)\in V'$
\begin{gather*}
g_{p,q}'(t)=f'(t|p|+(1-t)|q|)(|p|-|q|)+f'(t|q|+(1-t)|p|)(|q|-|p|)\\
g_{p,q}''(t)=f''(t|p|+(1-t)|q|)(|p|-|q|)^2+f''(t|q|+(1-t)|p|)(|q|-|p|)^2
\end{gather*}
Since $|p|>|q|$ and, by definition of admissibility, $f''(x)>0$
for every $x\in \mathbb R$, we conclude that $g_{p,q}''(t)>0$ for
every $t\in (0,1)$. Hence $g_{p,q}$ is strictly convex on $(0,1)$.
Moreover,
\[
g'_{p,q}(\frac{1}{2})=f'(\frac{|p|}{2}+\frac{|q|}{2})(|p|-|q|)+f'(\frac{|q|}{2}+\frac{|p|}{2})(|q|-|p|)=0.
\]
Since $g_{p,q}''>0$ on $(0,1)$, it follows that $g_{p,q}$ is
strictly convex on $(0,1)$ and achieves a unique minimum on
$(0,1)$ at $t=\frac{1}{2}$.

Since $0<P_f<\infty$ on $(0,1)$ and $\displaystyle
P_f=\sum_{(p,q)\in V'} g_{p,q}$, it also follows that $P_f$ is
strictly convex on $(0,1)$ and achieves a unique minimum on
$(0,1)$ at $t=\frac{1}{2}$.
\end{proof}

\section{Exploiting convexity}

In this section we combine the ideas of the previous two sections
and establish Theorem~\ref{C} from the introduction.

\begin{thm}\label{thm:conj2}
Let $k\ge 2$ and let $f:(0,\infty)\to (0,\infty)$ be monotone
decreasing function such that the following hold:
\begin{enumerate}
\item The function $f(x)$ is strictly convex on $(0,\infty)$.

\item
\[ \limsup_{x\to \infty} f(x)^{1/x}< \frac{1}{(2k-1)^k}.
\]
\end{enumerate}

Then there exists a convex neighborhood $U$ of $\mathcal L_\ast$
in $\Delta_k$ such that $0<P_f< C_f<\infty$ on $U$ and both $C_f$
and $P_f$ are strictly convex on $U$. In particular, $C_f\ne
const$ on $U$ and $P_f\ne const$ on $U$.
\end{thm}

\begin{proof}
By Theorem~\ref{thm:conj} there exists a
convex neighborhood $U$ of $\mathcal L_\ast$ in $\Delta_k$, such
that $0<P_f< C_f<\infty$ on $U$. We will prove that $P_f$ and
$C_f$ are strictly convex on $U$.

Let $D$ be the set of all $k$-tuples of integers $m=(m_1,\dots,
m_k)$ such that $m_i\ge 0$ for $i=1,\dots, k$ and $m_1+\dots
+m_k>0$. For each $m=(m_1,\dots, m_k)\in D$ let $Q_m$ be the set
of all $w\in \mathcal C_k$ such that $w$ involves exactly $m_i$
occurrences of $a_i^{\pm 1}$ for $i=1,\dots, k$ and let
$q_m:=\#(Q_m)$. Note that for every $w\in Q_m$, if $\mathcal
L=(x_1,\dots, x_k)\in \Delta_k$, then we have
\[
\ell_{\mathcal L}(w)=m_1x_1+\dots m_kx_k.
\]
Denote by $f_m:\Delta_k\to\mathbb R$ the function defined as
\[
f_m(x_1,\dots, x_k):=f(m_1x_1+\dots +m_kx_k), \quad (x_1,\dots,
x_k)\in \Delta_k.
\]
The function $f(x)$ is convex on $(0,\infty)$ and the function
$(x_1,\dots, x_k)\mapsto m_1x_1+\dots +m_kx_k$ is linear on
$\Delta_k$. Therefore $f_m$ is convex on $\Delta_k$.

Then for any $\mathcal L=(x_1,\dots, x_k)\in \Delta_k$ we have
\[
C_f(\mathcal L)=\sum_{m\in D}q_m f_m(\mathcal L).
\]
Since each $f_m$ is convex on $\Delta_k$, it follows that $C_f$ is
convex on $\Delta$. We claim that $C_f$ is strictly convex on $U$.
Let $D_1$ be the subset of $D$ consisting of all the $k$-tuples
having a single nonzero entry equal to $1$, that is, $D_1$ is the
union of the $k$ standard unit vectors in $\mathbb Z^k$. Let
$m_i=(0,\dots, 1, \dots 0)\in D_1$ where $1$ occurs in the $i$-th
position. Then $Q_{m_i}=\{[a_i], [a_i^{-1}]\}$ and $q_{m_i}=2$.
Also, $f_{m_i}(x_1,\dots, x_k)=f(x_i)$ for every $(x_1,\dots,
x_k)\in \Delta_k$.

Put $g:=f_{m_1}+\dots +f_{m_k}:\Delta_k\to\mathbb R$, so that
\[
g(x_1,\dots, x_k)=f(x_1)+\dots+f(x_k), \quad (x_1,\dots, x_k)\in
\Delta_k.
\]
It is easy to see that $g$ is strictly convex on $\Delta_k$ since
$f$ is strictly convex on $(0,\infty)$. We have:
\[
C_f=\sum_{m\in D}q_m f_m=2g+\sum_{m\in D-D_1}q_m f_m
\]
Since $C_f<\infty$ on a convex set $U$ and since $g$ is strictly
convex on $U$ and $\displaystyle\sum_{m\in D-D_1}q_m f_m$ is
convex on $U$, it follows that $C_f$ is strictly convex on $U$ as
claimed.

The proof that $P_f$ is strictly convex on $U$ is exactly the same
as for $C_f$ above. The only change that needs to be made is to
re-define $Q_m$ for each $m=(m_1,\dots, m_k)\in D$ as the set of
all $w\in \mathcal P_k$ such that $w$ involves exactly $m_i$
occurrences of $a_i^{\pm 1}$ for $i=1,\dots, k$.
\end{proof}

\begin{rem}
Let $\mathcal Z_k$ be the set of all \emph{root-free} conjugacy classes
$w\in \mathcal C_k$, that is, conjugacy classes of nontrivial elements
of $F_k$ that are not proper powers. It is not hard to show, similar to
Proposition~\ref{prop:ent}, that if $\mathcal L$ is a metric structure
on $W_k$ then
\[
h_\mathcal L=\widetilde h_{\mathcal L} 
\]
where
\[
\widetilde h_\mathcal L:=\lim_{R\to\infty}\frac{\log \#\{ w\in \mathcal Z_k:
\ell_{\mathcal
    L}(w)\le R \}}{R}.
\]
If one now re-defines the McShane function $C_f$ as $S_f$:
\[
S_f:\Delta_k\to (0, \infty], \qquad S_f(\mathcal L)=\sum_{w\in
\mathcal Z_k} f(\ell_\mathcal L(w))\quad \text{ where }\mathcal
L\in \Delta_k,
\]
then the proofs of the parts of Theorem~\ref{A} and Theorem~\ref{C}
dealing with $C_f$ go through verbatim for $S_f$.
\end{rem}

\end{document}